# Analytical Approach for Active Distribution Network Restoration Including Optimal Voltage Regulation

Hossein Sekhavatmanesh, *Student Member, IEEE*, Rachid Cherkaoui, *Senior Member, IEEE*

*Abstract*—The ever-increasing utilization of sensitive loads in the industrial, commercial and residential areas in distribution networks requires enhanced reliability and quality of supply. This can be achieved thanks to self-healing features of smart grids that already include the control technologies necessary for the restoration strategy in case of a fault. In this paper, an analytical and global optimization model is proposed for the restoration problem. A novel mathematical formulation is presented for the reconfiguration problem reducing the number of required binary variables while covering more practical scenarios compared to the existing models. The considered self-healing actions besides the network reconfiguration are the nodal load-rejection, the tap setting modification of voltage regulation devices (incl. OLTCs, SVR, and CBs), and the active/reactive power dispatch of DGs. The voltage dependency of loads is also considered. Thus, the proposed optimization problem determines the most efficient restoration plan minimizing the number of de-energized nodes with the minimum number of self-healing actions. The problem is formulated as a Mixed-Integer Second Order Cone Programming (MISOCP) and solved using the Gurobi solver via the MATLAB interface YALMIP. A real 83–node distribution network is used to test and verify the presented methodology.

## NOMENCLATURE

| | |
|---|---|
| ADN | Active Distribution Network |
| OLTC | On-Load Tap Changer |
| SVR | Step Voltage Regulator |
| CB | Capacitor Bank |
| SM | Smart Meter |
| DNO | Distribution Network Operator |
| DG | Distributed Generator |
| OPF | Optimal Power Flow |
| ENS | Energy Not Supplied |

### A. Parameters

| | |
|---|---|
| $w_{re}, w_{sw}, w_{op}, w_1, w_2$ | Weighting factors of the objective function terms (p.u.) |
| $D_i$ | Importance factor of the load at bus $i$ (p.u.) |
| $\lambda_{ij}$ | Weighting factor of the switch on line $ij$ depending on its operational priority |
| $\lambda_i$ | Weighting factor of the load breaker at bus $i$ depending on its operational priority |
| $P(/Q)^0_{i,t}$ | Active (/Reactive) power demand when voltage is 1 p.u. at bus $i$, at time $t$ (p.u.) |
| $(kp/kq)_i$ | Active (/Reactive) load voltage sensitivity for the load at bus $i$. |
| $r(/x)_{ij}$ | Resistance(/Reactance) of line $ij$ (p.u.) |
| $v^{max}/v^{min}$ | Maximum(/Minimum) limits of voltage magnitude (p.u) |
| $f_{ij}^{max}$ | Maximum current flow rating of line $ij$ (p.u) |
| $f_{ij}^{thr}$ | Current threshold level at line $ij$ beyond which the current deviation will be minimized (p.u) |
| $A_{i,p}$ | Indicator specifying if node $i$ is in zone $p$ (1/0) |
| $A_{ij,p}$ | Indicator specifying if line $ij$ is totally in zone $p$ (1/0) |
| $\sigma_i$ | Transformer ratio change at each tap step in the OLTC installed at bus $i$. |
| $\alpha 0_i$ | The pre-fault transformer ratio setting of OLTC/SVR at bus $i$ |
| $\Delta r0_i(/\Delta r0_{ij})$ | The pre-fault tap position of CB/SVR at bus $i$/line $ij$ |
| $\Delta Q_i$ | Reactive power change at each tap step in the CB installed at bus $i$ |
| $n_i(/n_{ij})$ | Number of tap steps at each positive/negative side of the OLTC/SVR/CB installed at bus $i$/line $ij$ |

### B. Variables

| | |
|---|---|
| $Y_{ij}$ | Binary decision variable indicating if the switched line $ij$ is energized or not (1/0) |
| $L_i$ | Binary decision variable indicating if the load at node $i$ is supplied or rejected (1/0) |
| $Z_{ij}$ | Continuous variable indicating if line $ij$ is oriented from node $i$ to node $j$ or not. |
| $\alpha_i$ | Continuous decision variable relaxing the post-fault transformer ratio setting of OLTC at bus $i$. |
| $\beta_{ij,t}$ | Continuous decision variable used in the linearization of nodal voltage equation for each SVR on line $ij$, at time $t$. |
| $b_k$ | The $k^{th}$-element of $\beta_{ij,t}$ binary representation |
| $\Delta r_i(/\Delta r_{ij})$ | Integer decision variable indicating the post-fault tap position of each CB/SVR at bus $i$/line $ij$. |
| $\delta r_k$ | The $k^{th}$-element of $\Delta r_{ij}$ binary representation |
| $\gamma_{ij,t}$ | Continuous variable approximating the tap position of each SVR on line $ij$, at time $t$ |
| $T_i(/T_{ij})$ | Continuous variable accounting for the tap changing of OLTC/CB/SVR at bus $i$/ line $ij$ |
| $E_p$ | Indicator of zone $p$, being energized or not (1/0) |
| $X_i$ | Indicator of node $i$, being energized or not (1/0) |
| $X_{ij}$ | Indicator of line $ij$, being energized or not (1/0) |
| $S_{ij}$ | Continuous variable indicating if sectionalizing switch on line $ij$ will be operated or not (1/0) |
| $B_i$ | Continuous variable indicating if load breaker at bus $i$ will be operated or not (1/0) |
| $F_{ij,t}$ | Square of current flow magnitude in line $ij$, at time $t$ (p.u) |
| $F_{ij,t}^*$ | Deviation of squared current flow magnitude in line $ij$, at time $t$ (p.u) |
| $P(/Q)_{i,t}^D$ | Active (/Reactive) power demand at bus $i$, at time $t$ (p.u.) |
| $P(/Q)_{i,t}^{cur}$ | Active/reactive load curtailment at bus $i$, at time $t$ (p.u.) |
| $p(/q)_{ij,t}$ | Active(/Reactive) power flow in line $ij$, starting from node $i$, at time $t$ (p.u) |
| $P(/Q)_{i,t}^{Sub}$ | Active(/Reactive) power injection from the substation node $i$, at time $t$ (p.u) |
| $Q_{i,t}^{inj}$ | Reactive power injected by a CB/DG at bus $i$, at time $t$ (p.u) |
| $P_{i,t}^{inj}$ | DG active power injection at bus $i$, at time $t$ (p.u) |
| $V_{i,t}$ | Square of voltage magnitude at bus $i$, at time $t$ (p.u) |

*C. Indices*

| | |
|---|---|
| $p$ | Index of zones |
| $i, j, k$ | Index of nodes |
| $ij$ | Index of branches |
| $t$ | Index of time |

*D. Sets*

| | |
|---|---|
| $N$ | Set of healthy nodes in the faulted feeder and available neighboring feeders |
| $N^*$ | Set of nodes in the off-outage area |
| $Z^*$ | Set of zones in the off-outage area |
| $W$ | Set of healthy lines in the faulted and available feeders |
| $W^*$ | Set of lines (plus tie-lines) in the off-outage area |
| $W^S = W_{ava}^S \cup W_{int}^S \cup W_{sec}^S$ | Set of lines (plus tie-lines) in the off-outage area equipped with switches |
| $W_{ava}^S$ | Set of lines hosting available tie-switches |
| $W_{int}^S$ | Set of lines hosting internal tie-switches |
| $W_{sec}^S$ | Set of lines hosting internal sectionalizing switches |
| $\Omega_{Sub}$ | Set of substation nodes in available feeders that are hosting an OLTC |
| $\Omega_{SVR}$ | Set of lines in the faulted and available feeders that are hosting an SVR |
| $\Omega_{CB}$ | Set of nodes in the faulted and available feeders that are hosting a CB |
| $\Omega_{DG}$ | Set of nodes in the faulted and available feeders that are hosting a DG |

## I. INTRODUCTION

In recent years, highly increasing penetration of DGs and voltage regulation devices such as OLTCs, SVRs, and CBs has brought significant challenges in ADN planning and operation [1, 2]. In parallel, smart grid concept evolved as a new operating paradigm for ADNs replacing conventional approaches that were deployed in passive distribution networks. One of the main features of smart grids is the self-healing capability. According to this feature, all the new monitoring, communication and control facilities pertaining to smart grids should be deployed resulting in an automatic and efficient restoration strategy [3].

Although distribution networks have a meshed structure, most of them are operated radially to ensure economic savings in terms of required number and size of protection devices. Following a failure in such a radial topology, once it is isolated, the area downstream to the fault place (*off-outage area*) remains unsupplied. This is the case, even when DGs exist in the off-outage area. Actually according to the IEEE standard 1547 [4], in case of a severe disturbance, every DG in the network should be automatically disconnected. The restoration service aims at supplying the off-outage area in the most efficient way. The objective is to pick up as many customers as possible with a minimum number of switching operations [5]. The new configuration of the network remains for a so-called *restorative period* until the faulted element is repaired. In the case of passive distribution networks, the restoration solution counts only on the switching operations. According to this strategy, the loads in the off-outage area are transferred to the healthy neighboring feeders through changing the status of normally-closed (sectionalizing) and normally-open (tie) switches.

Today, in ADNs, besides the switching operation, other self-healing actions can be deployed leading to a more efficient restoration solution [2]. Among these actions, one can cite the optimal set points of dispatchable DGs, OLTCs, AVRs, and CB. From the mathematical programming viewpoint, the restoration problem is a mixed-integer (due to the status of switches) and non-linear (due to the power flow constraints) optimization problem. In order to deal with the non-polynomial hardness of such an OPF-based optimization problem, some papers apply fuzzy algorithms [9–11] or meta-heuristic methods such as Genetic Algorithm (GA), Particle Swarm Optimization (PSO) and Tabu Search [9]. However, these methods are in general time-consuming and could fail to give a feasible solution in a reasonable time complying with online operation requirements. Therefore, heuristic approaches based on expert systems are proposed for the restoration problem to achieve a near-optimal solution in a short time [13–15]. These approaches use graph search methods suitable for specific network topologies and cannot be generalized for any network topology.

The authors in [14, 16] formulate the robust restoration problem as a mixed-integer programming problem. However, the active/reactive power losses are ignored in the linear formulation derived for the power flow constraints. Recently, some relaxations specific to radial distribution networks are proposed, making the OPF problem convex. These relaxations have been employed in the restoration problem, leading to an analytical optimization problem in the form of either a Mixed-Integer Quadratic Constraint Programming (MIQCP) [15] or a Mixed-Integer Second-Order Cone Programming (MISOCP) model [25, 26]. The convexity of the power flow formulation allows finding the global optimal solution up to the desired accuracy by using available commercial solvers.

However, the existing mathematical models for the restoration problem fail to represent partial restoration scenarios accurately. These scenarios are referred to the cases where the whole off-outage area cannot be restored while respecting all the electrical constraints (such as voltage and current limits). For dealing with this situation, the authors in [22, 23] propose to leave a part of the off-outage area without any supply which is called *isolated area*. In that work, in addition to the binary variables considered for the switch status, another set of binary variables is assumed to indicate if the corresponding nodes belong to the isolated area or not. Adding this auxiliary set of binary variables increases the complexity of the problem. Another approach for partial restoration scenarios is proposed in [20]. According to this method, the loads can be rejected in partial restoration scenarios but all the nodes in the off-outage area must be re-energized in any case. This could lead to increase uselessly the number of required switching operation.

In this paper, an analytical formulations in the form of a MISOCP is derived for the restoration problem, addressing the weaknesses of the existing methods in the modeling of partial restoration scenarios. The problem of computation burden is alleviated by decreasing the number of binary decision variables in the mathematical formulation together with reducing the border of the network part involved in the optimization problem. In this paper, the considered self-healing actions are the network reconfiguration and, if needed, load rejection, optimal active/reactive power set point of DGs and optimal tap settings of voltage control devices including OLTCs, SVRs, and CBs. The restoration strategy in this paper is to be implemented in a centralized way. However, it can also be implemented on each individual agent in a distribution network that is operated in a decentralized manner [21].

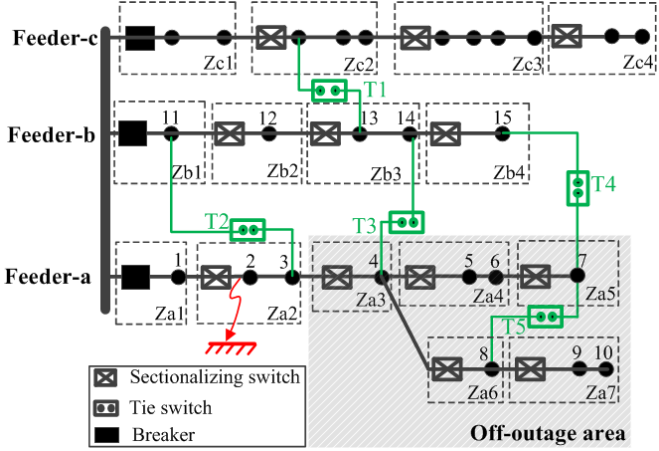

Fig. 1. A Simple schematic of a distribution network under a fault condition

## II. Problem Formulation

In this section, it is aimed to provide a convex analytical optimization model for the restoration problem while accounting for the operational limits using the full constraints of an AC-OPF problem. It is preferred for DNO to provide a network configuration for the restoration strategy closer to the pre-fault one. In this way, it will be easier to return to the topology of the normal state, once the fault is cleared. Therefore, in this work, only those healthy feeders that can be directly connected to the off-outage area through tie-switches are involved in the restoration process. These feeders and switches are called *available feeders* (Feeder-b in Fig. 1) and *available tie-switches* (T3 and T4 in Fig. 1), respectively. The tie-switches and sectionalizing switches with both ending nodes inside the off-outage area are referred as *internal tie-switches* (T5 in Fig. 1) and *internal sectionalizing switches*, respectively. It is assumed that the MV network that is under study is balanced and can be represented by a single-phase equivalent.

The number of discrete decision variables has the most weight in the computation burden of mixed-integer optimization problems. In the proposed model, the binary-or-integer decision variables are limited to 1) the load rejection status at each node in the off-outage area, 2) the tap position of SVRs and CBs in the off-outage area and available feeders, and 3) the energization status of available tie-switches, internal tie-switches, and internal sectionalizing switches.

The optimization problem is formulated as follows:

Minimize: $F^{obj} = w_{re}.F^{re} + w_{sw}.F^{sw} + w_{op}.F^{op}$ (1)

$$F^{re} = \sum_t \sum_{i \in N} D_i . (P_{i,t}^{cur} + Q_{i,t}^{cur})$$ (2)

$$F^{sw} = \sum_{(i,j) \in W_{ava}^S \cup W_{int}^S} Y_{ij}.\lambda_{ij} + \sum_{(i,j) \in W_{sec}^S} S_{ij}.\lambda_{ij} + \sum_{i \in N^*} B_i.\lambda_i$$ (3)

$$F^{Op} = w_1 . \sum_t \sum_{(i,j) \in W} F_{ij,t}^* + w_2 . \left( \sum_{i \in \Omega_{Sub} \cup \Omega_{CB}} T_i + \sum_{ij \in \Omega_{SVR}} T_{ij} \right)$$ (4)

Subject to:

$$0 \leq Z_{ij} \leq 1 \quad \forall (i,j) \in W^S \quad (5)$$

$$Z_{ij} = Y_{ij}, \ Z_{ji} = 0 \quad \forall (i,j) \in W_{ava}^S \quad (6)$$

$$Z_{ij} + Z_{ji} = Y_{ij} \quad \forall (i,j) \in W^S \quad (7)$$

$$\sum_{i \in N^*} \left( A_{i,p} . \sum_{j:(i,j) \in W^*} Z_{ji} \right) = E_p \leq 1 \quad \forall p \in Z^* \quad (8)$$

$$0 \leq Y_{ij} \leq M.F_{ij,t} \quad \forall (i,j) \in W^S, \forall t \quad (9)$$

$$\begin{cases} X_i = \sum_{p \in C^*} (A_{i,p}.E_p), & i \in N^* \\ X_i = 1, & i \in N \setminus N^* \end{cases} \quad (10)$$

$$\begin{cases} L_i \leq X_i, & i \in N^* \\ L_i = 1, & i \in N \setminus N^* \end{cases} \quad (11)$$

$$X_{ij} = \begin{cases} \sum_{p \in C^*} (A_{ij,p}.E_p), & (i,j) \in W^* \setminus W^S \\ Y_{ij}, & (i,j) \in W^S \\ 1, & (i,j) \in W \setminus W^* \end{cases} \quad (12)$$

*Switching operation constraints* (13)

*Tap setting constraints* (14)

*Load modelling constraints* (15)

*OPF constraints* (16)

The objective function given in (1) consists of reliability ($F^{re}$), switching ($F^{sw}$), and operational ($F^{op}$) terms, in decreasing order of priority. This hierarchical priority is enabled choosing appropriate values for $w_{re}$, $w_{sw}$, and $w_{op}$ accordingly using $\epsilon-constraint$ method [22]. The reliability objective term expressed in (2) tends to minimize the total active and reactive curtailed loads, while accounting for their importance factors. The constraints related to these curtailed powers are formulated in section II.D. The second objective term (3) corresponds to the switching criteria. It includes three sub-terms associated, respectively, to the operation of available and internal tie-switches, sectionalizing switches, and load breakers at each node in the off-outage area. With coefficient $\lambda_{ij}$ in (3), remotely controlled switches are prioritized for being operated over manually controlled ones. The operation of all these switches and breakers are determined using the constraints given in section II.A. The operational objective term is formulated in (4) through two sub-terms. The first one is the total deviation of squared magnitude of current flows beyond a certain threshold. This squared current flow deviation ($F_{ij,t}^*$) is formulated in the form of the following constraints [2].

$$F_{ij,t}^* \geq 0 \quad \forall (i,j) \in W, \forall t \quad (17)$$

$$F_{ij,t}^* \geq F_{ij,t} - (f_{ij}^{thr})^2 \quad \forall (i,j) \in W, \forall t \quad (18)$$

The second sub-term in (4) accounts for the tap changes in the voltage regulation devices. The corresponding constraints are formulated in details in section II.B. The relative importance of these two sub-terms is enabled using $w_1$ and $w_2$ depending on the operational policies of DNO during the restorative period. All the weighting factors are applied on the normalized terms of expressions defined in (2)-(4). In order to make the reliability and switching terms in the same scales (per unit), they are divided by their maximum possible values. The importance factors of loads and switches are taken into account to determine these maximum values.

Constraints (5)-(9) model the reconfiguration problem ensuring the radial topology of the network [23]. A binary decision variable $Y_{ij}$ is associated with each line in the off-outage area that is equipped with a switch indicating if it is energized or not. To each of these switched lines, is also assigned two continuous variables $Z_{ij}$ and $Z_{ji}$ indicating the orientation of the line with respect to the virtual source nodes. These nodes refer to the nodes outside the off-outage area that are connected to the available tie-switches (Ex. nodes 14 and 15 in Fig. 1). The presented constraints seem identical to those introduced in [16]. However, the formulation presented here accounts for the possible isolated areas in partial restoration scenarios. This claim is clarified in the proof of the following lemma.

**Lemma1:** in the feasible solution space, Z variables must take only binary values, and, $Y_{ij} = 1$ shows that line $ij$ is an edge of a connected spanning tree.

**Proof:** Let $W^y = \{(i,j) \in W^S | Y_{ij} = 1\}$ is the set of energized switched lines corresponding to a specific solution for $Y$, according to (5)-(9). Consider bus 15 that is a virtual source node in case of the fault shown in Fig. 1. If line (15,7) belongs to $W^y$, then $Z_{15-7} = 1$ according to (6). It means that we have one entering flow to zone[1] $Za5$ hosting node 7. Thus, by (8), there should not be any other entering flow[2] to zone $Za5$. It means that $Z_{6-7}$ must be zero. Therefore, according to (7) if line $(7,6) \in W^y$, $Z_{7-6} = 1$. In other words, zone $Za4$ that is hosting node 6 is not feeding zone $Za5$ but is supplied from zone $Za5$. The same analysis made for zone $Za5$ can be applied for zone $Za4$ and other zones in the paths originating from bus 15. From such a recursive analysis it can be concluded that each line through $W^y$ is oriented such that its Z variable toward the source node is 0, and the one outgoing from the source node is 1. Therefore, all the branches of $W^y$ that are in a path connected to bus 15 takes the form of a connected tree without forming any loop. If one part of the network does not contain any path connected to a source node, it means that its buses are not receiving any supply. Consequently, the current flow variables must be 0 ($F_{ij,t} = 0$). This enforces $Y_{ij}$ variables for the lines in such an isolated area to be 0, according to (9), M being a large multiplier. We conclude that $W^y$ is composed of spanning trees, each originated from one source node. Thus, the set of variables $Y$ results in a radial configuration of the network with the possibility of having isolated areas in the network. ∎

According to (10), a node in the off-outage area is energized ($X_i = 1$) if it is inside an energized zone ($E_p = 1$). For such an energized node in the off-outage area, a decision is made in (11) with binary variable $L_i$, indicating if its load will be restored or rejected (1/0). Outside the off-outage area, the loads are kept connected all the time ($L_i = 1$). According to (12), energized switched lines are identified directly with variable $Y$, while the other lines will be energized if their hosting zones are energized.

### A. Switching operation constraints

As expressed in (4) normally open switches including available and internal tie-switches are operated when they are energized ($Y_{ij} = 1$). However, for sectionalizing switches, energization status does not necessarily imply that they should be operated or not. For example, in case of a partial restoration scenario, sectionalizing switches that are entirely in the isolated part are de-energized, but do not need to be opened. As formulated in (13.a), the sectionalizing switch on line $ij$ must be operated (opened) only if it is de-energized ($Y_{ij} = 0$) and at least one of its ending buses are energized ($X_i = 1$ or $X_j = 1$). According to (13.a), $S_{ij}$ for other sectionalizing switches can get a value between 0 and 1. However, since minimization of $S_{ij}$ is included in the objective function, $S_{ij}$ for these switches will be zero, meaning that they should not be opened. As given in (13.b), the load breaker at each MV/LV substation node should be operated only if the node is energized ($X_i = 1$) and the load is to be rejected ($L_i = 0$). Otherwise, the breaker will not be operated ($B_i = 0$) according to the problem solution, since its value is to be minimized as a term in the objective function (4).

$$\begin{cases} (1 - Y_{ij}) + X_i - 1 \leq S_{ij} \leq 1 \\ (1 - Y_{ij}) + X_j - 1 \leq S_{ij} \leq 1 \qquad \forall (i,j) \in W_{sec}^S \quad (13.a) \\ 0 \leq S_{ij} \end{cases}$$

$$(1 - L_i) + X_i - 1 \leq B_i \leq 1 \qquad \forall i \in N^* \quad (13.b)$$

### B. Tap setting constraints

Three types of voltage regulation devices in distribution networks are studied in this work, namely OLTCs at the HV/MV substation nodes ($\Omega_{Sub}$), SVRs on certain branches ($\Omega_{SVR}$), and CBs at some nodes ($\Omega_{CB}$). In this part, a linear modeling is proposed for these devices considering a limited number of discrete decision variables. The aim is to find their optimal tap positions regarding the voltage control objectives while minimizing the number of required tap changes for the restorative period. It is assumed that the tap setting of each voltage regulator remains unchanged for the whole restorative period.

#### 1) On-Load Tap Changing Transformers[3]

The voltage at the HV side of OLTC is assumed to be constant (1 p.u.), while the voltage at the MV side is regarded as the slack voltage in the OPF constraints formulated in section 0. The OLTC transformer ratio changes a little between two consecutive tap positions (e.g., $\sigma = 0.625\%$ as reported in [20]). Therefore, in order to relax the problem complexity, we represent the step-wise values of tap ratio $\Delta r \in \{-n.\sigma, \ldots, -\sigma, 0, \sigma, \ldots, n.\sigma\}$ with a continuous variable $\alpha \in [-n.\sigma, n.\sigma]$. Finally, the optimal value obtained for $\alpha$ will be rounded to the closest $\Delta r$ value in the original discrete set. With this assumption, the constraints related to each OLTC installed at a HV/MV substation are formulated as follows:

$$V_{i,t} = 1.(1 + \alpha_i)^2 \approx 1 + 2\alpha_i \qquad \forall i \in \Omega_{Sub}, \forall t \quad (14.a)$$

$$-n_i.\sigma_i \leq \alpha_i \leq n_i.\sigma_i \qquad \forall i \in \Omega_{Sub} \quad (14.b)$$

---

[1] In this paper, a zone is referring to each segment of the feeder that is surrounded by two or more sectionalizing switches (Fig. 1).

[2] It is meant directional flow with regard to the graph theory

[3] OLTC and SVR are modeled as ideal transformers in series with virtual lines representing their impedances.

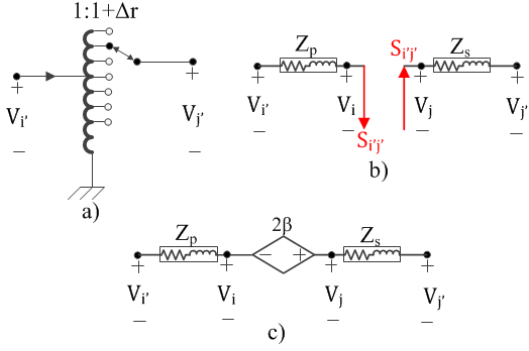

Fig. 2. Modeling of SVR in the restoration problem. a) SVR schematic b) standard equivalent model c) linearized equivalent model

$$\begin{cases} T_i \geq \dfrac{\alpha_i - \alpha 0_i}{\sigma_i} - 1 \\ T_i \geq -1 - \dfrac{\alpha_i - \alpha 0_i}{\sigma_i} \\ T_i \geq 0 \end{cases} \quad \forall i \in \Omega_{Sub} \quad (14.c)$$

Regarding the small positive values of $\alpha_i$, squared voltage relationship in (14.a) is linearized using binomial approximation. As expressed in the second term of (5), it is aimed for each OLTC to minimize the change of its tap position with respect to the tap position prior the fault occurrence ($\alpha 0_i$). The tap position will be changed only if the deviation of relaxed transformer ratio ($\alpha_i$) from the initial one ($\alpha 0_i$) is higher than the ratio step size ($\sigma_i$). This condition is linearized in (14.c) defining a set of auxiliary continuous variables ($T_i$) showing quantitatively the extent of tap change.

*2) Step Voltage Regulator*

Fig. 2(b) shows the equivalent model of the SVR shown in Fig. 2(a), where $Z_p$ and $Z_s$ denote the impedances at the primary and secondary sides, respectively. Considering this equivalent model, the optimal tap setting constraints are formulated in the following:

$$V_{j,t} = V_{i,t} \cdot (1 + \Delta r_{ij})^2 \quad \forall (i,j) \in \Omega_{SVR}, \forall t \quad (14.d)$$
$$\approx V_{i,t} \cdot (1 + 2\Delta r_{ij}) = V_{i,t} + 2\beta_{ij,t}$$

$$\begin{cases} T_{ij} \geq \Delta r_{ij} - \Delta r 0_{ij} \\ T_{ij} \geq \Delta r 0_{ij} - \Delta r_{ij} \\ T_{ij} \geq 0 \end{cases} \quad \forall (i,j) \in \Omega_{SVR} \quad (14.e)$$

Similar to OLTCs, the granularity of tap positions in SVRs allows using the binomial approximation. However, since the voltage at the starting node cannot be assumed 1 p.u. as for OLTCS, tap ratio setting cannot be represented by continuous variables. An auxiliary set of continuous variables ($\beta_{ij,t}$) is introduced to linearize the product of $\Delta r_{ij} \cdot V_{i,t}$. For this aim, the integer variable $\Delta r_i \in \{-n, n\}$ is represented by 2n+1 binary variables denoted by $\delta r_k$, as given in (14.f) and (14.g). Constraints (14.h) enforce continuous variables $b_k$ to take the value of $\delta r_k \cdot V$, M being a large multiplier. Accordingly, the value of $\Delta r \cdot V$ denoted by $\beta$ is retrieved from weighted summation of $b_k$ variables (14.i).

$$\Delta r = \sum_k k \cdot \delta r_k \quad : \quad k \in \{-n, \dots, -1, 1, \dots, n\} \quad (14.f)$$

$$\sum_k \delta r_k = \delta r_0 \quad : \quad k \in \{-n, \dots, -1, 1, \dots, n\} \quad (14.g)$$

$$\begin{cases} b_k \leq M \cdot \delta r_k \\ b_k \leq V \\ b_k \geq V - M \cdot (1 - \delta r_k) \\ b_k \geq 0 \end{cases} \quad \begin{array}{l} \forall k \\ \in \{-n, \dots, -1, 1, \dots, n\} \end{array} \quad (14.h)$$

$$\beta = \sum_k k \cdot b_k \quad : \quad k \in \{-n, \dots, -1, 1, \dots, n\} \quad (14.i)$$

Adding binary variables ($\delta r_k$) does not increase the computation cost too much, since the above-mentioned constraints are in the form of the Special Ordered Set-1 (SOS1) constraints [24]. With these types of constraints, the number of nodes to be searched in the underlying branch and bound process will be reduced significantly. Finally, the model shown in Fig. 2(c) will be the linearized model of each SVR that can be incorporated into the OPF constraints presented in section II.D. In this model, the tap changing amount will be modeled as expressed in (14.e).

*3) Capacitor Banks*

Regarding the CBs, the amount of reactive power change at each step ($\Delta Q_i$) is not small. Therefore, unlike OLTCs and SVRs, the step changes of CBs should be modeled with integer variables ($t_i$). Therefore, the model of CBs participating in the restoration process is as follows:

$$Q_{i,t}^{inj} = \Delta Q_i \cdot \Delta r_i \quad \forall i \in \Omega_{CB}, \forall t \quad (14.j)$$

$$\Delta r_i \in \{0,1,2, \dots, n_i\} \quad \forall i \in \Omega_{CB} \quad (14.k)$$

$$\begin{cases} T_i \geq \Delta r_i - \Delta r 0_i \\ T_i \geq \Delta r 0_i - \Delta r_i \\ T_i \geq 0 \end{cases} \quad \forall i \in \Omega_{CB} \quad (14.l)$$

Likewise OLTCs and SVRs, the objective related to the operation of CBs is to minimize the tap changing which is modeled in (14.l) using an auxiliary set of continuous variables ($T_i$).

*C. Load modeling constraints*

In order to find the optimal set point for the regulation devices during the restorative period, it is important to consider the voltage dependency of load consumption. In this paper, the exponential model is used for the active/reactive loads as given in (15.a) and (15.b).

$$P_{i,t}^D = P_{i,t,s}^0 \left(\dfrac{v_{i,t}}{v_0}\right)^{kp_i} \quad \forall i \in N, \forall t \quad (15.a)$$

$$Q_{i,t}^D = Q_{i,t}^0 \left(\dfrac{v_{i,t}}{v_0}\right)^{kq_i} \quad \forall i \in N, \forall t \quad (15.b)$$

This model is linearized assuming that $v_0 = 1\, p.u.$ and $V_{i,t,s}$ is close to $1\, p.u.$, so that the binomial approximation can be applied according to the following:

$$P_{i,t}^D = P_{i,t}^0 (1 + (V_{i,t} - 1))^{kp_i/2}$$
$$\approx P_{i,t}^0 (1 + \dfrac{kp_i}{2}(V_{i,t} - 1)) \quad \forall i \in N, \forall t \quad (15.c)$$

$$Q_{i,t}^D = Q_{i,t}^0 (1 + (V_{i,t} - 1))^{kq_i/2}$$
$$\approx Q_{i,t,s}^0 (1 + \dfrac{kq_i}{2}(V_{i,t} - 1)) \quad \forall i \in N, \forall t, s \quad (15.d)$$

*D. OPF Constraints*

In the following, the relaxed formulation of AC-OPF is presented, incorporating the network reconfiguration problem,

optimal tap setting of voltage regulators and modeling of load-voltage dependencies. The aim is to re-dispatch the active/reactive power set points of DGs during the restorative period, while respecting all the security constraints in the reconfigured network.

$$0 \leq P_{i,t}^D - P_{i,t}^{cur} \leq M.L_i$$
$$0 \leq Q_{i,t}^D - Q_{i,t}^{cur} \leq M.L_i \quad \forall i \in N, \forall t \quad (16.a)$$

$$0 \leq P_{i,t}^{cur} \leq M.(1-L_i)$$
$$0 \leq Q_{i,t}^{cur} \leq M.(1-L_i) \quad \forall i \in N, \forall t \quad (16.b)$$

$$0 \leq F_{ij,t} \leq X_{ij}.f_{ij}^{max2} \quad \forall (i,j) \in W, \forall t \quad (16.c)$$

$$-M.X_{ij} \leq p_{ij,t} \leq M.X_{ij}$$
$$-M.X_{ij} \leq q_{ij,t} \leq M.X_{ij} \quad \forall (i,j) \in W, \forall t \quad (16.d)$$

$$v_i^{min2}.X_i \leq V_{i,t} \leq v^{max2}.X_i \quad \forall i \in N, \forall t \quad (16.e)$$

$$-M.(1-X_{ij}) \leq V_{i,t} - V_{j,t} - 2(r_{ij}.p_{ij,t} + x_{ij}.q_{ij,t})$$
$$\leq M.(1-X_{ij}) \quad \forall (i,j) \in W \backslash \Omega_{SVR}, \forall t \quad (16.f)$$

$$p_{ij,t} = \left( \sum_{\substack{i^* \neq i \\ (i^*,j) \in W}} p_{ji^*,t} \right) + r_{ij}.F_{ij,t} + P_{j,t}^D - P_{j,t}^{cur} - P_{j,t}^{Sub}$$
$$- \left( P_{j,t}^{inj} : j \in \Omega_{DG} \right) \quad \forall t, (i,j) \in W \quad (16.g)$$

$$q_{ij,t} = \left( \sum_{\substack{i^* \neq i \\ (i^*,j) \in W}} q_{ji^*,t} \right) + x_{ij}.F_{ij,t} + Q_{j,t}^D - Q_{j,t}^{cur} - Q_{j,t}^{Sub}$$
$$- (Q_{j,t}^{inj} : j \in \Omega_{CB} \cup \Omega_{DG}) \quad \forall t, (i,j) \in W \quad (16.h)$$

$$\left\| \begin{matrix} 2p_{ij,t} \\ 2q_{ij,t} \\ F_{ij,t} - V_{i,t} \end{matrix} \right\|_2 \leq F_{ij,t} + V_{i,t} \quad \forall (i,j) \in W, \forall t \quad (16.i)$$

$$0 \leq P_{i,t}^{inj} \leq P_{i,max}^{inj}$$
$$Q_{i,min}^{inj} \leq Q_{i,t}^{inj} \leq Q_{i,max}^{inj} \quad \forall i \in K_4, \forall t \quad (16.j)$$

$$\left\| \begin{matrix} P_{i,t}^{inj} \\ Q_{i,t}^{inj} \end{matrix} \right\|_2 \leq S_{i,max}^{inj} \quad \forall i \in K_4, \forall t \quad (16.k)$$

In all the OPF constraints, $M$ is a fixed as a large multiplier. Constraints (16.a) and (16.b) enforce active/reactive load curtailment variables ($P_{i,t,s}^{cur}/Q_{i,t,s}^{cur}$) being equal to the active/reactive load powers for de-energized loads ($L_i = 0$) and zero for energized ones ($L_i = 1$). Constraints (16.c) and (16.d) enforce, respectively, current and active/reactive power flows to be zero for each de-energized line. The squared voltage magnitude is set to zero, in (16.e), for de-energized nodes ($X_i = 0$), and kept within the feasible region for energized nodes. The nodal voltage constraint in [16] is revised in (16.f) to exclude unrestored lines. Note that for the links representing the ideal SVR (ex. line $ij$ in Fig. 2C) (14.d) is applied instead of (16.f).

Constraints (16.g) and (16.h), respectively, concern with the active-/reactive power balances at the end buses of each line, taking into account the active/reactive load curtailments and power injections by CBs and DGs. Constraint (16.i) is the relaxed version of the current flow equation in each line proposed in [16]. Constraint (16.j) limits the active/reactive power injection of each dispatchable DG, while its apparent power is limited by the cone constraint given in (16.k).

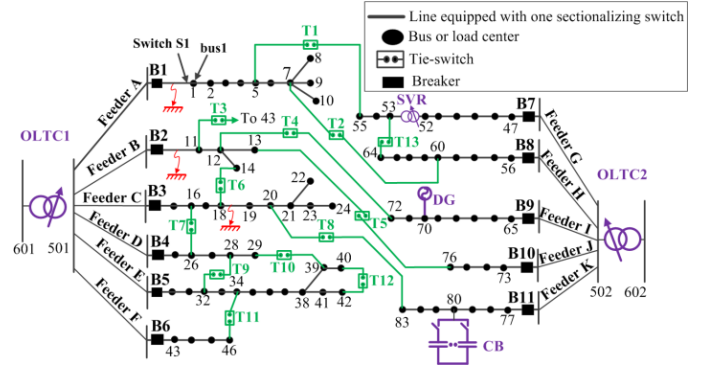

Fig. 3. Test distribution network with 4 feeders and 70 buses [25]

## III. SIMULATION AND RESULT

The proposed mathematical solution model of the restoration problem is tested on a 11.4 kV distribution network shown in Fig. 3. This test system is based on a practical distribution network in Taiwan with 2 substations, 4 feeders, 83 nodes, and 96 branches (incl. tie-branches) [26]. Each line is assumed to be equipped with a sectionalizing switch at the receiving node. The detailed nodal and branch data is given in [26]. The base power and energy are chosen as 1 MVA and 1 Mwh, respectively. The network is updated adding two OLTCs at the both substations, one SVR on line 52-53, one CB at bus 80, and one dispatchable DG at bus 70. The set point of these devices are determined by DNO a day ahead, which might be subject to further tuning for the restoration strategy. The current threshold for each line ($f_{ij}^{thr}$ in (4)) is assumed to be 50% of its current capacity limit. According to ANSI C84.1 standard, the minimum and maximum voltage magnitude limits for the restorative period are set, respectively, to 0.917 and 1.050 p.u. [27].

The majority of nodes along feeders A, B and E in Fig. 3 are assigned, respectively, to rural, light and high industrial load types. Most of the nodes in feeders C, F, H, J have load characteristics close to residential loads and most of the nodes along feeders D, G, I, and K supply commercial loads. There is also 10-30% of public light loads considered at each node. The hourly profile for each of these load patterns is given in [25]. Tie-switch T4 and sectionalizing switch S6 are remotely controlled and the others are all manually controlled. The algorithm is implemented on a PC with an Intel(R) Xeon(R) CPU and 6 GB RAM; and solved in Matlab/Yalmip environment, using Gurobi solver.

The proposed restoration strategy is applied on the test system shown in Fig. 3 according to five simulation scenarios. The optimal solution of the restoration strategy in each simulation scenario is reported in Table I. In all the cases, the restorative period is assumed from 8:00 A.M. to 10:00 P.M., during which most of the loads experience their peak values. It is assumed that the default priority factors of loads at buses 3 and 5 are 10 and for all the other loads are 1. The load voltage sensitivity coefficients ($K_p$ and $K_q$ in (15.c) and (15.d)) are assumed the same for all the loads. The default values of $K_p$ and $K_q$ are equal to 0.6 and 0.3, respectively.

Table I. Optimal restoration results in case of different fault scenarios

| Simulation Scenarios | Fault location | Sequence of switching | ENS (p.u.) | Total current deviation (p.u.) | Optimal setting of voltage regulation devices | Computation time (sec) |
|---|---|---|---|---|---|---|
| 1 | 1-2 | - Open S6 and S7<br>- Close T1 and T2 | 13.22 | 7.77 | - | 3.78 |
| 2a | 1-2 | - Open S3 and S7<br>- Close T1 and T2 | 1.34 | 69.98 | $\frac{v_{502}}{v_{602}} = 1.025$ | 2.54 |
| 2b | 1-2 | - Open S3 and S7<br>- Open load breaker 5<br>- Close T1 and T2 | 4.56 | 46.73 | $\frac{v_{502}}{v_{602}} = 1.019$ | 2.68 |
| 3 | 1-2 | - Open S7<br>- Close T1 and T2 | 0 | 82.59 | $\frac{v_{502}}{v_{602}} = 1.025$<br>$\frac{v_{53}}{v_{52}} = 1.0125$ | 4.59 |
| 4a | 18-19 | - Open load breaker 20<br>- Close T8 | 2.02 | 27.60 | $\frac{v_{502}}{v_{602}} = 1.025$ | 0.90 |
| 4b | 18-19 | - Close T8 | 0 | 10.97 | $\frac{v_{502}}{v_{602}} = 1.025$<br>$Q_{80}^{inj} = 0.6\ p.u.$ | 0.30 |
| 5a | 11-12 | - Close T5 | 0 | 5.80 | $\frac{v_{502}}{v_{602}} = 1.0125$ | 2.13 |
| 5b | 11-12 | - Close T4 | 0 | 42.49 | $\frac{v_{502}}{v_{602}} = 1.025$ | 2.51 |

*A. Scenario1: a partial restoration scenario*

Scenario 1 is defined where a fault occurs on line 1-2, as shown in Fig. 3. In this scenario, it is assumed that none of the voltage regulation devices are controllable, fixing their tap positions at zero. As it can be noticed from the results, node 6 is isolated from the other parts of the unsupplied area that are restored through tie-switches T1 and T2. This type of maneuver is missing in the restoration strategies proposed in the literature. Actually, according to those strategies, instead of opening the remotely-controlled switch S7, the load breaker at node 6 must be disconnected. This leads to increase uselessly the deployment time of the restoration strategy. It is also worth to mention that in this scenario, instead of load at node 6, loads at buses 3, 4 and 5 could have been rejected, decreasing the total ENS value from 13.22 p.u. to 10.77 p.u. However, regarding the higher load priority factors at nodes 3 and 5, they are preferred to be restored. In addition, instead of opening S6, the load breaker at bus 6 could have been opened. However, operating S6, which is remotely controlled helps to accelerate the restoration solution deployment.

*B. Scenario2: effects of OLTC optimal tap setting*

In scenario 2a, the restoration strategy in case of the fault on line 1-2 is studied while incorporating the tap control of OLTCs in the two substations. These OLTCs enable ±5% voltage regulating range, in 8 steps ($\sigma = 0.625\%, n = 4$). The results of this scenario are shown in Table I. With the optimal setting of OLTC-2, instead of isolating node6 as in scenario 1, node 2 is isolated, decreasing the ENS by one order of magnitude. The optimal voltage ratio of OLTC-2 is obtained 3.84 ($\frac{\alpha_i - \alpha 0_i}{\sigma_i}$ in (14.c)), corresponding to a tap position equal to +4 (the closest discrete value). As it can be seen, the computation time remains very small, thanks to the continuous representation of OLTC tap positions.

In order to see the effect of load voltage sensitivity coefficients ($K_p$ and $K_q$ in (15.c) and (15.d)), the simulation scenario 2.b is studied. For this aim, $K_p$ and $K_q$ are increased from their default values. They are equal to 2 for all the loads along feeders H, I, J, and K. According to the results shown in Table I, the optimal restoration strategy in scenario 2.b is to isolate bus 2 and reject load at bus 5 leading to a larger value of ENS compared to scenario 2a. As it can be seen, the optimal tap position of OLTC is changed from 4 in scenario 2.a to 3 (2.75 which is rounded to 3) in scenario 2.b. Comparing the results of these two scenarios highlights a trade-off between keeping the voltage profile at the extremities above the minimum limit and limiting the increase of the load consumption level.

*C. Scenario3: effects of SVR optimal tap setting*

In scenario3, it is assumed that the fault occurs on line 1-2 and the SVR on line 52-53 is controllable besides the OLTC-2. This SVR enables ±5% voltage regulation range in 8 steps ($\sigma = 0.625\%, n = 4$). According to the results shown in Table I, incorporating optimal tap setting of SVR in the restoration strategy helps to restore all the nodes. The optimal tap position of SVR is obtained as +2. As it can be seen, the computation time still remains reasonable.

*D. Scenario4: effects of CB optimal tap setting*

In this part, the effects of incorporating optimal tap setting of CB on the restoration strategy is studied. First of all, scenario 4.a is considered, where the restoration problem in case of a fault on line 18-19 is studied in the same condition as in scenario 2a. As reported in Table I, the load at node 20 cannot be restored in these conditions. The effect of optimal CB tap setting is illustrated in scenario 4b. In this scenario, the fault on line 18-19 is considered again while incorporating the CB at node 80, with 1.2MVAR injection capacity in 4 steps. As it can be seen in Table I, putting the tap of this CB on +2 enables to restore the whole off-outage area. The detailed results of scenario 4b is reported in Table I.

*E. Scenario5: effects of optimal DG incorporation*

In this part, the restoration problem is studied in case of the fault on line 11-12. First, scenario 5a is defined such as among all the voltage control devices, only OLTC-2 is under control. In this condition, as shown in Table I, the optimal solution is to pick up the whole unsupplied area through tie-switch T5. T5 is selected because of the third priority term in the objective function. Operating T6 instead of T5 leads to a total current deviation value equal to 46.73 p.u., which is larger than 5.80 p.u. This scenario is to be compared with scenario 5b, where the DG at node 70 is incorporated into the restoration problem. The active and apparent power capacities of DG at node 70 are 2.8 and 3.0 p.u., respectively. As it can be seen in Table I, in scenario 5b, with the help of DG power injection (mainly at t=11:00 with reactive power) switch T4 is able to restore the whole unsupplied area. Although operating switch T5 leads to a very lower total current deviation in the network, switch T4 is chosen, since it is remotely-controllable. Note that the switching term is given more priority compared to the technical term. The DG active/reactive power set-point and the resulting apparent power dispatch is shown in Fig. 4. As it can be seen, the active and apparent power limits of DG are well respected.

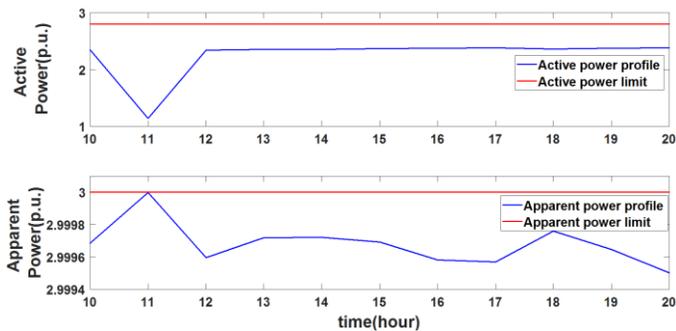

Fig. 4. Optimal set-points for DG active and apparent power during the restorative period

## IV. Conclusion

This paper presents a MISOCP optimization model for the restoration problem in an ADN. A novel mathematical model is proposed for the reconfiguration problem. This model covers among others partial restoration scenarios that were missing in the literature. The convex relaxation methods for the OPF problem is applied in the proposed restoration problem ensuring that the electrical constraints are all respected during the whole restorative period. Besides the network reconfiguration, the optimal setting of DGs and voltage regulation devices are incorporated into the problem as other self-healing measures. The proposed restoration model was successfully tested on a real 83-bus distribution network in five simulation scenarios. The results confirm the contributions claimed by the authors demonstrating the effectiveness of the proposed strategy in terms of quality and computation performance.